\documentclass[12pt]{article}
\usepackage{amsmath, amssymb}
\usepackage[dvipdfmx]{graphicx}
\usepackage[dvipdfmx]{color}
\usepackage{enumerate}
\usepackage[numbers,sort]{natbib}
\usepackage{comment}
\usepackage{url} 

\usepackage{amsthm}
\theoremstyle{definition}
\newtheorem{theorem}{Theorem}

\newtheorem{proof of lemma}{proof of Lemma}
\usepackage{here}

\newcommand{\blind}{1}

\begin{document}

\def\spacingset#1{\renewcommand{\baselinestretch}%
{#1}\small\normalsize} \spacingset{1}

\if1\blind
{
 \title{\bf Asymptotic Analysis of the Bayesian Likelihood Ratio for Testing Homogeneity in Normal Mixture Models}
 \author{Natsuki Kariya, and Sumio Watanabe\hspace{.2cm}\\
 Department of Mathematical and Computing Science, \\ Tokyo Institute of Technology\\
 }
 \maketitle
} \fi

\if0\blind
{
 \bigskip
 \bigskip
 \bigskip
 \begin{center}
 {\LARGE\bf Asymptotic Analysis of Singular Likelihood Ratio of Normal Mixture  by Bayesian Learning Theory for Testing Homogeneity}
\end{center}
 \medskip
} \fi

\bigskip

\begin{abstract}
When we use the normal mixture model, the optimal number of the components describing the data should be determined. Testing homogeneity is good for this purpose;
however, to construct its theory is challenging, since the test statistic does not converge to the $\chi^{2}$ distribution even asymptotically. 
The reason for such asymptotic behavior is that 
the parameter set describing the null hypothesis (N.H.) 
contains singularities in the space of the alternative hypothesis (A.H.). 
Recently, a $\it{Bayesian}$ theory for singular models was developed, and it has elucidated various problems of statistical inference. 
However, its application to hypothesis tests for singular models has been limited. 

In this paper, we introduce a scaling technique that greatly simplifies the derivation and study
testing of homogeneity for the first time the basis of Bayesian theory.  
We derive the asymptotic distributions of the marginal likelihood ratios in three cases:
 (1) only the mixture ratio is a variable in the A.H. ;
 (2) the mixture ratio and the mean of the mixed distribution are variables;
 And  (3) the mixture ratio,  the mean, and the variance of the mixed distribution are variables.;
In all cases, the results are complex, but can be described as functions of random variables obeying normal distributions. A testing scheme based on them was constructed, and their validity was confirmed through numerical experiments. 
\end{abstract}

\noindent%
\vfill

\newpage
\spacingset{1.45} 
\section{Introduction}
\label{sec:intro}

Normal mixtures have been widely used for analyzing various problems, such as pattern recognition, clustering analysis, and anomaly detection, since they were first applied to Pearson's biology research
 in the 19th century \citep{Pearson1894}. This remains one of the most important models in statistics, both in theory and in practice \citep{Mclachlan2000}.

When a normal mixture model is employed, the optimal number of components for describing the data has to be determined.  Testing homogeneity is a well-known approach for this purpose, which is a hypothesis test
 to determine whether the data are described by a single normal distribution or a mixture distribution.
In a normal mixture, the correspondence between a parameter and a probability density function is not one-to-one, and the Fisher information matrix of the statistical model that represents alternative hypotheses becomes singular at the parameter of the 
null hypothesis. As a result, the log likelihood ratio of the test of homogeneity for the normal mixture model does not converge to a $\chi^{2}$ distributions, unlike the regular models \citep{Hartigan1985}\citep{Liu2003}\citep{Garel2001}. 

Therefore, it is necessary to study the testing of homogeneity in normal mixture models, not only out of theoretical interests but also for practical applications. Various methods have been proposed; for example, the modified likelihood ratio test, a method that adds a regularizing term \citep{Chen2001}\citep{Chen2004}, an EM algorithm for calculating the modified likelihood ratio \citep{Chen2009}\citep{Chen2012}, and the D test \citep{Charnigo2004}.  (for a recent review on this topic, see for example \citep{Chauveau2017}). 
However, little research  exists on treating the problem from the \textit{Bayesian} perspective.

On the other hand, a theoretical foundation for singular statistical models has been constructed within the framework of Bayesian statistics in recent years\citep{Watanabe2018}. One of the achievement of this theoretical study is WAIC, a new information criterion that can be applied to singular models\citep{Watanabe2010}. However, most of the results have been on the problem of statistical inference,  whereas the problem of hypothesis testing remains insufficiently studied. 

In this paper, we study the test of homogeneity of normal mixture models based on the framework of a Bayesian hypothesis test, for the first time.
We derive the asymptotic distribution of the test statistic, i.e.,the marginal likelihood ratio, in three cases: 
 (1) only the mixture ratio is a variable; (2) the mixture ratio and the mean of the mixed distribution in the A.H. are variables; (3) the mixture ratio,  the mean, and the variance of the mixed distribution in the A. H. are variables. In all cases, the marginal likelihood ratios converge to certain random variables, which are different from the well-known $\chi^{2}$ distribution as an effect of the singularities in the model. 
The validity and efficiency of the derived theory are shown numerically. 

The paper is organized as follows. In Section 2, we review the framework of the Bayesian hypothesis test
and show that the marginal likelihood ratio gives the most powerful test.  Our main results are presented in Sections 3 to Section 5. We derive the asymptotic distributions of the marginal likelihood ratio analytically for 
three cases. The results of the numerical experiment for validation are also presented.  In Section 6, we summarize our results and give a conclusion.

\section{Framework of Bayesian Hypothesis Test}

\label{sec:Bayesian}

In this section, we briefly review the framework of the Bayesian hypothesis test. 

Let $\left\{X^{n}\right\} =(X_1,X_2,...,X_n)$ be a sample which is generated independently and identically
from a probability distribution. 
We consider a statistical model of a normal mixture, 
\begin{equation}
p(x|w) = (1-a) \mathcal{N}(0,1^2) + a\mathcal{N}(b,\frac{1}{c}),
\end{equation}
where $w=(a,b,c)$, $0\leq a \leq 1$, $b\in\mathbb{R}$, and $c > 0$. 
Here ${\cal N}(b,\sigma^2)$ denotes a normal distribution with the average $b$ and 
 variance $\sigma^{2}$. 

In a Bayesian hypothesis test, the null and alternative hypotheses are set as
\begin{eqnarray*}
\mbox{N.H.} & : & w_0\sim \varphi_0(w),\;\;\;X_i\sim p(x|w_0), \\
\mbox{A.H.} & : & w_0\sim \varphi_1(w),\;\;\;X_i\sim p(x|w_0),
\end{eqnarray*}
where $X \sim p(x)$ means that a random variable $X$ is generated from a probability density function
 $p(x)$. 
In the case of the testing homogeneity, the null hypothesis is set as,
\begin{eqnarray*}
\mbox{N.H.} & : & \varphi_0(w) = \delta(a)\delta(b)\delta(c-1).
\end{eqnarray*}
For a given statistic $T(X^n)$ and a real value $t$, a hypothesis test $T$  is defined by a determining procedure, 
\begin{eqnarray*}
T(X^n)\leq t &\Longrightarrow & \mbox{N.H.}, \\
T(X^n)>t & \Longrightarrow & \mbox{A.H.}.
\end{eqnarray*}
The level and power of this hypothesis test are defined by the probability that the A.H. is chosen on the assumption that  the N.H. and A.H. generate $X^n$ respectively. 
\begin{eqnarray*}
\mbox{Level}(T,t)&=& \mbox{Probability}(\mbox{A.H.}|\mbox{N.H.} ), \\
\mbox{Power}(T,t)&=& \mbox{Probability}(\mbox{A.H.}|\mbox{A.H.} ). 
\end{eqnarray*}
For given two hypothesis tests $T$ and $U$,  $T$ is said to be more powerful than $U$ if 
and only if 
\[
\mbox{Level}(T,t)=\mbox{Level}(U,u)
\Longrightarrow 
\mbox{Power}(T,t)\geq\mbox{Power}(U,u)
\]
holds for an arbitrary set $(t,u)\in\mathbb{R}^2$. 
A test $T$ is said to be most powerful  if it is more powerful than any other test. 

In regard to Bayesian hypothesis test, it was proved that a test using the marginal likelihood ratio as a test statistic is the most powerful test, where the marginal likelihood ratio is defined as
\begin{equation}
L(X^{n}) = \frac{\displaystyle \int \varphi_{1}(w) \prod_{i=1}^np(X_{i}|w) dw}
{\displaystyle \int \varphi_{0}(w) \prod_{i=1}^np(X_{i}|w) dw}. 
\end{equation}
where  $X^n=\{X_i\}$ is an i.i.d. sample. 
Therefore, the probability 
distribution of the marginal likelihood ratio $L(X^n)$ is necessary for constructing
 the most powerful test in the Bayesian framework.

Note that this test is different from the Bayesian model selection using the marginal likelihood ratio, 
because this test is NOT defined by choosing the A.H. when $L(X^n)\geq 1$. 
The most powerful test of this paper is defined by choosing the A.H. when $L(X^n)\geq t$ for 
$t$ which makes $P(A.H.|N.H.)$ be a given level. 

In this paper, we mathematically derive the asymptotic probability distributions of $L(X^{n})$ in the 
following three cases for A.H., 
\begin{enumerate}
\item $\varphi_{1}(a,b,c) = U_{a}(0,1) \delta(b-\beta)\delta(c-1),$ 
\item $\varphi_{1}(a,b,c) = U_{a}(0,1)  U_{b}(0,B)\delta(c-1),$
\item $\varphi_{1}(a,b,c) = U_{a}(0,1)   U_{bc}(D),$
\end{enumerate} 
where $U_a(0,1)$ is the uniform distribution of $a$ on the interval $(0,1)$, $U_b(0,B)$ is the uniform distribution of $b$ on $(0,B)$, and $U_{bc}(D)$ is the uniform distribution of a 
set $D$ in $(b,1/c)$ space. 

The proofs use the following notation. 
For a given sample $X^n$, two random variables $\xi_n$ and $\eta_n$ are
defined by
\begin{eqnarray}
\xi_n&=& \frac{1}{\sqrt{n}}\sum_{i=1}^n X_i, 
\label{eq:xin}
\\
\eta_n&=&\frac{1}{\sqrt{2n}}\sum_{i=1}^n (X_i^2-1).
\label{eq:etan}
\end{eqnarray}
If $X^n$ is an i.i.d. sample generated from the N.H., then 
both $\xi_n$ and $\eta_n$ converge to $\mathcal{N}(0,1)$ in 
distribution and they are asymptotically independent.

\section{Case 1: the case only the mixture ratio is unknown}

\subsection{Asymptotic distribution of the test statistic}

Let us consider the case 1,  that is, the mean of the mixed distribution of the A.H. is fixed
 and only the mixture ratio is the variable. 
This case is quite simple and we can readily derive the asymptotic distribution of the marginal likelihood ratio analytically. However, even in such a case  the marginal likelihood ratio shows non-conventional behavior (different from the ordinary $\chi^{2}$ distribution), as the support of the prior 
of  the A.H. approaches to the singularity in the parameter space. 

Therefore, this case can be regarded as a minimal model with which to study the effect of the singularity on the behavior of the marginal likelihood ratio. Hence we will study it as a first step towards analyzing more practical situations in the following sections. 

We consider the case that the A.H. is near the N.H. in terms of the Kullback-Leibler divergence. In this situation, it is not easy to discriminate the alternative hypothesis from the null one. This is a typical situation in which a hypothesis test is needed.  

A similar situation occurs in the context of the Bayesian $\it{inference}$, where the true distribution generating the sample is slightly deviates from the singularity of the model on the order of $O(n^{-1/2})$ is studied\citep{Watanabe2003} . Here, it was shown that  the singularity greatly affects the behavior of the generalization error, even when the parameter set that represents the true model does not definitely match the singularity. 

Although our problem is not an inference but a hypothesis test,
we expected that a similar structure exists.
 We will see that this is true, and that the $n^{-1/2}$ scaling works as well. This is because the scaling is determined from the order of the Kullback-Leibler divergence between the
 A.H. and the singularity (N.H.).

Applying the scaling mentioned above, we can derive the asymptotic distribution of the marginal likelihood ratio as follows.
\begin{theorem}\label{theorem:1}
Assume that the N.H. and A.H. are given as 
\begin{eqnarray*}
\mbox{N.H.} & : & \varphi_0(w) = \delta(a)\delta(b)\delta(c-1), \\
\mbox{A.H.} & : & \varphi_1(w) = U_a(0,1)\delta(b-\beta)\delta(c-1),
\end{eqnarray*}
where $\beta = \beta_{0}\times n^{-\frac{1}{2}}$ and $\beta_{0}$ is a nonzero constant. 
If $\{X_i\}$ is independently and identically generated from the N.H., 
the convergence in probability,
\[
L(X^n)-L_{\infty}(\xi_n)\rightarrow 0
\]
holds for $n\rightarrow \infty$, where
\begin{eqnarray}
L_{\infty}(\xi_n) &= & \frac{\sqrt{2\pi}}{2\beta_{0}}\left[{\rm erf}
\left(\frac{\beta_0-\xi_n}{\sqrt{2}}\right)
+{ \rm erf}\left(\frac{\xi_n}{\sqrt{2}}\right)\right]
 \exp(\frac{\xi_n^{2}}{2}).
\end{eqnarray}
Here $\xi_n$ is a random variable defined in eq.(\ref{eq:xin}) and
${\rm erf}(x)$ is the error function,
\begin{equation*}
{\rm erf}(x) \equiv \frac{2}{\sqrt{\pi}}\int_{0}^{x} e^{-t^{2}}dt.
\end{equation*}
\end{theorem}
Remark. Assume that $\xi$ is a random variable whose probability
distribution is $\mathcal{N}(0,1)$. By Theorem \ref{theorem:1} 
and the convegence in distribution $\xi_n\rightarrow\xi$, 
the convergence in distribution $L(X^n)\rightarrow L_{\infty}(\xi)$ holds.
Since $L_{\infty}(\xi)$ can be rewriten as 
\[
L_{\infty}(\xi)=\int_0^1da\; \exp
\left(-\frac{\beta_0^2a^2}{2}+\beta_0\xi a\right),
\]
it is an increasing function of $\xi$. Therefore, 
we can determine the rejection region by using $\xi$. 

\begin{proof}:
The integral with respect to  $(b,c)$ is easily performed and the prior of $a$ in the A.H. is a uniform distribution
on $(0,1)$; it follows that 
\begin{eqnarray}
L(X^{n}) & = & \int_0^1 \exp(H(a))\;da,
\end{eqnarray}
where $H(a)$ is defined as, 
\begin{eqnarray}
H(a)&\equiv&
\sum_{i=1}^n\log{\frac{p(X_{i}|a,\beta,1)}{p(X_{i}|0,0,1)}}\\
&=&
\sum_{i=1}^n\log{\left\{(1-a) + a\exp{\left(\beta X_{i}-\frac{\beta^{2}}{2}\right)}\right\}}.
\label{eq:logL}
\end{eqnarray}
Under the N.H., 
from a well-known result in extreme statistics, the order of the maximum of $X_{i}$ is
\begin{equation*}
X_{M}\equiv \max{\left\{X_{i}\right\}}  = O_{p}\left(\sqrt{2\log{n}}\right).
\end{equation*}
This results in 
\begin{equation}
\beta X_{i} - \frac{\beta^{2}}{2} \leq \beta X_{M}-\frac{\beta^{2}}{2} 
= O_{p}\left(\sqrt{\frac{\log{n}}{n}}\right).
\end{equation}
Let $\alpha$ be a constant which satisfies $1<\alpha$. Then
\begin{equation*}
\beta X_{i} - \frac{\beta^{2}}{2} \sim o_{p}\left(\sqrt{\frac{\left(\log{n}\right)^{\alpha}}{n}}\right).
\end{equation*}
Hence
\begin{eqnarray*}
\exp\left(\beta X_{i} - \frac{\beta^{2}}{2}\right)&=&
1+\left(\beta X_{i}-\frac{\beta^{2}}{2}\right)+\frac{1}{2!}\left(\beta X_{i}-\frac{\beta^{2}}{2}\right)^{2}
\\&&+\frac{1}{3!}\left(\beta X_{i}-\frac{\beta^{2}}{2}\right)^{3}\times e^{C_{0}},
\end{eqnarray*}
where $C_0$ is a random variable that satisfies
\[
|C_0|\leq \left|\beta X_{i} - \frac{\beta^{2}}{2}\right|.
\]
Then,
\[
 \begin{cases}
 0 \leq C_{0} \leq \beta X_{i}-\frac{\beta^{2}}{2} & 
 (\mbox{ if } \beta X_{i}-\frac{\beta^{2}}{2} \geq 0), \\
 \beta X_{i}-\frac{\beta^{2}}{2} \leq C_{0} \leq 0 & (\mbox{otherwise}).
 \end{cases}
\]
Therefore,
\begin{equation*}
\frac{1}{3!}\left(\beta X_{i}-\frac{\beta^{2}}{2}\right)^{3}e^{C_{0}} \sim o_{p}\left(\frac{\left(\log{n}\right)^{3\alpha/2}}{n^{\frac{3}{2}}}\right). 
\end{equation*}
It follows that 
\begin{eqnarray*}
H(a)&=& \sum_{i=1}^n\log{
\left[
1+a\left\{
 \left(
 \beta X_{i}-\frac{\beta^{2}}{2}
 \right)+\frac{1}{2}
 \left(
 \beta X_{i}-\frac{\beta^{2}}{2}
 \right)^{2}
 \right\}+ o_{p}(\frac{1}{n})
\right]
} \\
&=& \sum_{i=1}^n\log{
\left[
1+a
 \beta X_{i}-\frac{a\beta^{2}}{2}
+\frac{a \beta^2 X_{i}^2}{2}
 + o_{p}(\frac{1}{n})
\right].
} 
\end{eqnarray*}
Then, by applying a Taylor expansion $\log(1+\epsilon)=\epsilon-\epsilon^2/2+O(\epsilon^3)$ to this equation, we obtain
\begin{equation}
H(a) = \sum_{i=1}^n \left[a\beta X_{i} - \frac{1}{2}a\beta^{2}+\frac{1}{2}a\beta^{2}X_{i}^{2}-\frac{1}{2}a^{2}\beta^{2}X_{i}^{2}
\right]+o_p(1). 
\end{equation}
Let us use the following notations,
\begin{eqnarray*}
\gamma &\equiv& \frac{\sum_{i}\beta X_{i} + \frac{1}{2}\sum_{i}\left(\beta X_{i}\right)^{2}-\frac{1}{2}\beta^{2}}{\frac{1}{2}\sum_{i}\left(\beta X_{i}\right)^{2}}, \\
\delta &\equiv& \frac{1}{2}\sum_{i}\left(\beta X_{i}\right)^{2}.
\end{eqnarray*}

Accordingly, $H(a)$ can be written as 
\begin{eqnarray*}
H(a)&=&-\delta a^2+\gamma\delta a
\\
&=&
-\delta (a-\gamma/2)^2+\delta\gamma^2/4.
\end{eqnarray*}
It follows that 
\begin{eqnarray*}
L(X^{n}) &=& \int_{0}^{1} da \exp{\left[-\delta(a-\frac{1}{2}\gamma)^{2}\right]} \times\exp{\left[\frac{1}{4}\times \gamma^{2}\delta \right]}\\
& = &\frac{\sqrt{\pi}}{2\sqrt{\delta}}\left[\rm{erf}\left(\frac{\gamma \sqrt{\delta}}{2}\right) + \rm{erf}\left(\sqrt{\delta}(1-\frac{\gamma}{2})\right)\right]\times \exp{\left[\frac{1}{4}\times \gamma^{2}\delta \right]},
\end{eqnarray*}
where $\rm{erf}(x)$ is the error function defined by 
\begin{equation*}
{\rm erf}(x) = \frac{2}{\sqrt{\pi}}\int_{0}^{x} e^{-t^{2}}dt.
\end{equation*}
As  $n$ tends to infinity, $\delta$ converges in probability as
\begin{equation}
\delta \rightarrow \frac{1}{2}\beta_{0}^{2},
\end{equation}
by which $\gamma$ satisfies 
\begin{equation}
\gamma = 2 \xi_n+o_p(1),
\end{equation}
which completes the theorem. 
\end{proof}

A remarkable feature of this theorem is that $ L (X ^ {n}) $ does not explicitly depend on the sample size $n$. The reason is that,  in the current setting, the distance between two centers of the clusters is $O(n^{-1/2}) $, and as the sample size $n$ increases, the posterior distribution becomes localized around the true parameter. But at the same time, the fluctuation around the true parameter induced by the randomness of the sample is of the same magnitude as the speed that the posterior distribution approaches the true parameters as the sample size increases. As a result of this, these two effects cancel and $ L (X ^ {n}) $ does not explicitly depend on the sample size $n$.

\subsection{Numerical evaluation of the level}
Here, we numerically derive the rejection region and the level based on the results above. 
From the definition, the level of the test is given as the probability that $ L (X ^ {n}) $ exceeds a certain threshold value $a$.

To see its behavior, we numerically calculated the level by generating 1,000 random samples from a standard normal distribution $\mathcal {N} (0,1^{2}) $ and calculated $ L (X^{n})$ by using them according to each sample. Then, we evaluated the level as the portion of the $ L (X^{n})$ that exceeded the threshold. 

\begin{figure}
\begin{center}
\includegraphics[width=10cm]{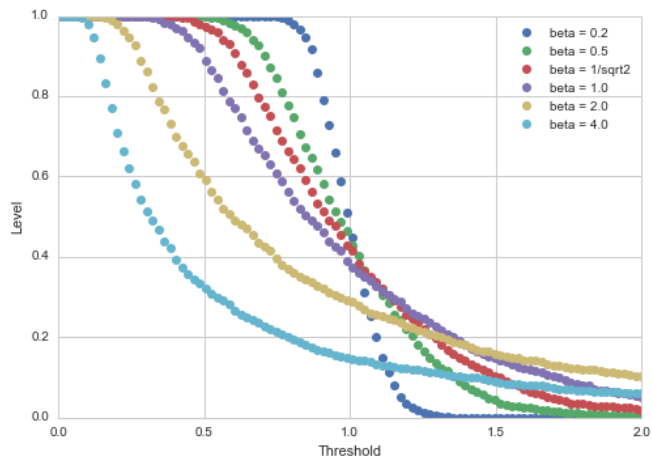}
\caption{Level calculated from the asymptotic distribution of $L$ as a function of the threshold for each $\beta$}
\label{thres-level}
\end{center}
\end{figure}

Figure \ref {thres-level} shows the plot of the level as a function of the threshold for each $\beta$.

The level drops rapidly when the threshold exceeds a certain value. This tendency is especially clear in small $\beta$ cases. This can be understood as follows.

The situation we consider is a ''delicate``, one in which it is not easy to discriminate between the N.H. and the A.H. is not clear. If we choose a small $a$ threshold, the level is large, but supposing that we choose larger and larger values, the level sharply decreases, because the N.H. and the A.H. become similar in the situation, and the probability that the value of marginal likelihood ratio is large is expected to be very low.

\begin{figure}[h]
\begin{center}
\includegraphics[width=10cm]{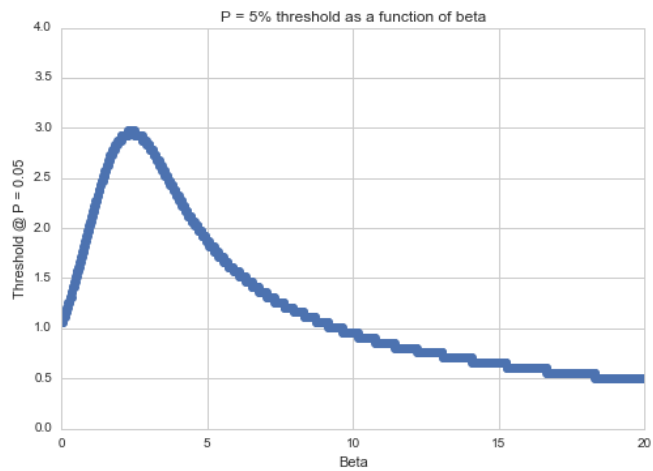}
\caption{Relation between the threshold that gives a 5\% level and $\beta$}
\label{beta-thres}
\end{center}
\end{figure}

Figure \ref{beta-thres} shows the results of our numerical calculation of the threshold that gives the 5\% level as a function of $\beta$. 
From the asymptotic distribution obtained above,  it can be seen that $ L \rightarrow \frac {1} {2} $ when $\beta $ is sufficiently large, and $ L \rightarrow 1 $ is within the limit of $\beta \rightarrow 0 $. 

To confirm the validity of the above result obtained above, we numerically evaluated the value of $ L (X ^ {n}) $ for a finite number of samples. We used a set of samples generated from the N.H. distribution. Here, we set the sample size as $n = 50 $ or $n = 100 $, and generated $10000$ sets of samples from the null hypothesis. In each case, the level was calculated from them as the ratio of the number of sets $X ^ {n}$ for which $ L (X ^ {n})$ falls into the critical region to the total number of sets. 

The result is shown in the Table 1, where we can see that the numerically calculated levels match those derived from the asymptote in both $n = 50$ and $n = 100 $ cases. 
Therefore,  it can be concluded that the asymptotic distribution we derived above is valid.

\begin{table}[h]
\begin{center}

\caption{Level calculated from the samples generated from the null hypothesis}

\vspace{5mm}
 
 \begin{tabular}{|c|c||r|r|} \hline
 \multicolumn{2}{|c||}{parameters}
 & \multicolumn{2}{|c|}{level} \\\hline 
 $\beta$ & threshold of 5\% & n = 50 & n = 100 \\\hline
 0.5 & 1.464 & 5.07\%& 4.76\% \\\hline
 1 & 2.02 & 5.56\%& 5.15\% \\\hline
 1.5 & 2.475 & 5.13\%& 4.92\% \\\hline
 2 & 2.929 & 4.88\%& 5.03\% \\\hline
 \end{tabular}
 \label{tab:label}
 \end{center}
\end{table}

\section{Case 2: both the mixture ratio and the mean of the mixed distribution are unknown}
\label{sec:AsympUnif}

In this section,  we proceed towards a bit more general case.  Here, the alternative hypothesis is a normal mixture whose mixture ratio and the mean are both variables. We assume that we have no prior knowledge about the A.H., and that the prior is a uniformly distribution. 

\subsection{Asymptotic distribution of the test statistic}

We will prove the following theorem on the asymptotic distribution of the marginal likelihood ratio $L$.
 
\begin{theorem}\label{theorem:2} 
Assume that N.H. and A.H. are given by
\begin{eqnarray*}
\mbox{N.H.} & : & \varphi_0(w) = \delta(a)\delta(b)\delta(c-1), \\
\mbox{A.H.} & : & \varphi_1(w) = U_a(0,1)U_b(0,B)\delta(c-1),
\end{eqnarray*}
respectively, where
$
B = B_{0}\times n^{-\frac{1}{2}}.
$
Then, if $X^n$ is an i.i.d. sample generated from N.H., the convergence 
in probability
\[
\
L(X^n) - L_{\infty}(\xi_n)\rightarrow 0
\]
holds as $n\rightarrow\infty$, where 
\begin{eqnarray}
L_{\infty}(\xi_n) &=& \frac{1}{B_{0}}\int_{0}^{B_{0}^{2}}\frac{1}{2\sqrt{t}}\log{\left(\frac{B_{0}^{2}}{t}\right)}e^{-t/2}\cosh{\left(\xi_n\sqrt{t}\right)}dt,
\label{eq:Linfty2}
\end{eqnarray}
and $\xi_n$ is a random variable defined in eq.(\ref{eq:xin}). 
\end{theorem}
Remark. Assume that $\xi$ is a random variable whose probability
distribution is $\mathcal{N}(0,1)$. By Theorem \ref{theorem:2},
the convergence in distribution
\[
L(X^n)\rightarrow L_{\infty}(\xi)
\]
holds. Hence the asymptotic rejection region of the most powerful test
can be found by using $L_{\infty}(\xi)$. 
\begin{proof}: 
The marginal likelihood ratio can be written as  
\begin{eqnarray*}
L(X^{n}) &=& \int_{0}^{1} da \int_{0}^{B} \frac{db}{B}\;
\exp(H(a,b)), 
\end{eqnarray*}
where
\begin{eqnarray*}
H(a,b)=\sum_{i}\log{
\left\{
(1-a) + a\exp{\left(b X_{i}-\frac{1}{2}b^{2}\right)}
\right\}.
}
\end{eqnarray*}
From the condition $b\in[0,B_0/\sqrt{n}]$, the $H(a,b)$ can be approximated in the same way as in
 the proof of Theorem 1 as,
\begin{equation}
H(a,b)=\sum_{i} \left[ab X_{i} - 
\frac{1}{2}ab^{2}+\frac{1}{2}ab^{2}X_{i}^{2}-\frac{1}{2}a^{2}b^{2}X_{i}^{2}\right]
+o_p(1).
\end{equation}
Hence,
\begin{equation}
H(a,b)= -\frac{n}{2}a^{2}b^{2} + \sum_{i}\left[abX_{i} + \frac{1}{2}ab^{2}\left(X_{i}^{2}-1\right)-\frac{1}{2}a^{2}b^{2}\left(X_{i}^{2}-1\right)\right]
+o_p(1).
\end{equation}
Under the N.H., 
by using the definitions, eqs.(\ref{eq:xin}) and (\ref{eq:etan}),  we have
\begin{eqnarray*}
H(a,b)&=& -\frac{n}{2}a^{2}b^{2} + \sqrt{n}
\left(
ab\xi_n+ \frac{1}{2}ab^{2}\eta_n-\frac{1}{2}a^{2}b^{2}\eta_n
\right)
+o_p(1) \\
&=& -\frac{n}{2}a^{2}b^{2} + \sqrt{n}
ab\xi_n
+o_p(1).
\end{eqnarray*}
Using the notation $L=L(X^n)$ for simplicity, it follows that 
\begin{eqnarray}
L &=& \int_{0}^{1} da \int_{0}^{B_{0}/\sqrt{n}}\frac{\sqrt{n}\;db}{B_{0}}  \;
\exp(-\frac{n}{2}a^{2}b^{2} + \sqrt{n}ab\xi_n+o_p(1))\nonumber 
\\
&=& \int_{0}^{1} da \int_{0}^{B_{0}}\frac{db}{B_{0}} \;
\exp(-\frac{1}{2}a^{2}b^{2} + ab\xi_n+o_p(1)). \nonumber 
\end{eqnarray}
Hence, the convergence in probability $L(X^n)-L_\infty(\xi_n)\rightarrow 0$ holds, where
\begin{eqnarray}
L_{\infty}(\xi_n)= \int_{0}^{1} da \int_{0}^{B_{0}}\frac{db}{B_{0}} \;
\exp(-\frac{1}{2}a^{2}b^{2} + ab\xi_n) \nonumber 
\end{eqnarray}
By using $b=t/a$, we have
\begin{eqnarray}
L_{\infty}(\xi_n)&= & \int_{0}^{1} da \int_{0}^{a B_{0}}
\frac{dt}{a B_{0}} \;
\exp(-\frac{1}{2}t^2 + t\xi_n) \nonumber 
\\
&= &  \int_{0}^{B_{0}} \frac{dt}{B_{0}} \int_{t/B_0}^{1} 
\frac{da}{a} \;
\exp(-\frac{1}{2}t^2 + t\xi) \nonumber 
\\
&= &\frac{1}{B_{0}} \int_{0}^{B_{0}} dt\left( \log(B_0) -\log t\right)
\exp(-\frac{1}{2}t^2 + t\xi). \nonumber 
\\
&= &\frac{1}{2B_{0}} \int_{0}^{B_{0}} dt\left( \log(B_0) -\log t\right)
\exp(-\frac{1}{2}t^2)\cosh( t\xi). \nonumber 
\end{eqnarray}
Then eq.(\ref {eq:Linfty2}) is obtained by replacing the integration of $t$ by $\sqrt{t}$. 
\end{proof}

Similar to the previous example, we can see that the asymptotic behavior of the test statistics $L$ does not explicitly depend on $n$.

We should also notice that the stochastic behavior of $L$ is determined only by that of the random variable $\xi$. Clearly, $L$ increases monotonously as the absolute value of $\xi$ increases, and $\cosh{\left(\xi\sqrt{t}\right)}$ is an even function with respect to $\xi$, hence, we can determine the critical region in the same way as is done in a two-sided hypothesis test of $\xi$.

For example, under the null hypothesis, the random variable $\xi$ obeys the standard normal distribution, and the 5 \% critical region is given as $|\xi| > 1.96$.

As a result of this, the 5 \% critical region of the test statistics $L$ is given as follows, 

\begin{equation}
|L|> \frac{1}{B_{0}}\int_{0}^{B_{0}^{2}}\frac{1}{2\sqrt{t}}\log{\left(\frac{B_{0}^{2}}{t}\right)}e^{-t/2}\cosh{\left(1.96 \sqrt{t}\right)}dt.
\end{equation}

For example, if we choose $B_{0} = 1$, the 5\% critical region of $L$ is given as

\begin{equation*}
|L| > 2.298
\end{equation*} 

We numerically validated the effectiveness of the analytically derived distribution of $L$ when the sample size is finite.  

First, we prepared the 10000 sets of the $n$ samples, where $n$ means the sample size and we set $n$ as $n=50$ or $n=100$.
We calculated the $L$ by substituting the $\xi$ in the asymptote with $\frac{1}{\sqrt{n}}X_{i}$. Here, we fixed $B_{0}$ as 1.
In each case, the level was calculated from them as the ratio of the number of sets $X ^ {n}$ for which $ L (X ^ {n})$ falls into the critical region to the total number of sets.  The levels were compared with those calculated from the level calculated from the asymptote of $L$. 

Table 2 shows the result. It shows the asymptote we derived in the previous section works well even in the finite $n$ cases.

\begin{table}[h]
\begin{center}

\caption{The level calculated from the samples generated from the null hypothesis}

\vspace{5mm}
 
 \begin{tabular}{|c|c|r|r|} \hline
 level & 10\% & 5\% & 1\% \\\hline
 rejection region  $L>r$ & r=2.171 & r=2.298 & r=2.646 \\\hline
numerically calculated level(n=50) & 9.75\% & 5.22\%& 1.04\% \\\hline
numerically calculated level(n=100) & 9.91\% & 4.73\%& 0.97\% \\\hline
numerically calculated level(n=200) & 9.74\% & 4.84\%& 0.99\% \\\hline
 \end{tabular}
 \label{tab:label}
 \end{center}
\end{table}

\subsection{Comments on the comparison with hypothesis test using Bayes factor}

Let us comment on the comparison of our results with those obtained by another well-known method of Bayesian hypothesis testing, i.e.,using the Bayes factor.  

As we derive the asymptote of the marginal likelihood ratio $L$, we can readily calculate the log marginal likelihood ratio $F$.

\begin{equation*}
F = -\log{L}
\end{equation*}

The log marginal likelihood ratio $F$ , which is also called the logarithm of the Bayes factor, can be used as a tools for hypothesis testing. 
The procedure is very simple and effective,  and it is used in various situations. 

The procedure is as follows. When the value of $F$ calculated from the data becomes negative, we choose the alternative hypothesis, and if otherwise, we choose the null hypothesis otherwise.

For the present problem, we can consider two ways of hypothesis testing with the result we derived.  One is based on the stochastic behavior of $L$, and the other is based on the $F$.
Both use the same quantity $L$, but we will see below that the former may work more effective in the  ``delicate`` situation.

\begin{figure}
\begin{center}
\includegraphics[width=10cm]{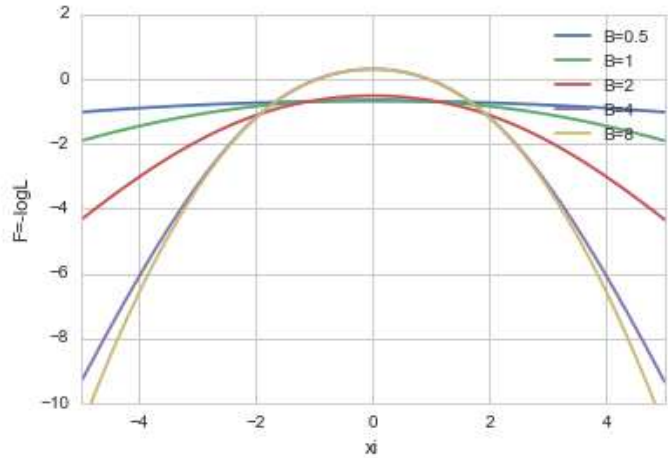}
\caption{The log marginal likelihood ratio $F$ as a function of the random variable $\xi$ for several values of $B_{0}$.}
\label{xi-logL}
\end{center}
\end{figure}

Figure \ref{xi-logL} shows the behavior of $F$ as a function of $\xi$.

Interestingly, when $B_{0}$ is small, the value of $F$ is always negative, regardless of any $\xi$, while in the large  $B_{0}$ case, $F$ becomes positive in the small $\xi$ region.
This can be understood as follows. When the two centers of the mixture distribution are so close as the distance between them is $\mathcal{O}(n^{{-1/2}})$, the overlap of the distribution of the null hypothesis and the distribution of the alternative hypothesis is large, and the sign of Bayes factor can become negative for any $\xi$.

In other words, when the two hypotheses are difficult to distinguish, the hypothesis test using the Bayes factor may choose the alternative hypothesis for any data, and it does not work well. On the other hand, the likelihood ratio test based on the stochastic behavior of $L$ is expected to work in such delicate cases.

\section{Case 3: the case the mixture ratio, the mean of the distribution mixed, and the variance are unknown}

\subsection{Asymptotic distribution of test statistic}

Here, we discuss a more practical case in which the variance of the A.H. is also a variable. 
That is, we consider the following probabilistic model,
\begin{equation}
p(x|a,b,c) = \left(1-a\right)\mathcal{N}(0,1) + a \mathcal{N}(b, \frac{1}{c}).
\end{equation}
We set the N.H. and the A.H. as
\begin{eqnarray*}
\mbox{N.H.} &:&  \varphi_0(a,b,c) = \delta(a) \delta(b)\delta(c-1), \\
\mbox{A.H.} &:&   \varphi_1(a,b,c)  = U_{a}(0,1) U(b,c),
\end{eqnarray*} 
where $U_a(0,1)$ is a uniform distribution on the interval $(0,1)$, and 
$U(b,c)$ is a uniform distribution on an ellipsoid in the $(b,c)$ plane such as,
\[
D=\left\{(b,c)\;;\;b^{2} + \frac{(c-1)^{2}}{2} \le \frac{R_0^2}{n}\right\},
\]
where  $R_{0}$ is a constant. The area of $D$ is $\sqrt{2}\pi R_0^2/n$. 

\begin{theorem}\label{theorem:3}
When the sample size $n\rightarrow\infty$, 
convergence in probability
$L(X^n)- L_{\infty}(\Xi_n)\rightarrow 0$ holds, where
\begin{eqnarray*}
L_{\infty}(\Xi_n)=\frac{1}{2R_0^2}\int_0^{R_0^2} 
\left(
\frac{R_0}{\sqrt{t}}-1
\right)
\exp\left(-\frac{t}{2}\right)
I_0(\sqrt{t}\Xi_n)\;dt.
\end{eqnarray*}
Here,  $I_{0}(t)$ is the modified Bessel function, 
\[
I_0(t)=\frac{1}{\pi}\int_0^\pi\cosh(t\sin\theta)d\theta, 
\]
which is monotone increasing in $t>0$
and $\Xi_n$ is a random variable defined by
\begin{equation*}
\Xi_n =\sqrt{\xi_n^2+\eta_n^2},
\end{equation*} 
where $\xi_n$ and $\eta_n$ are defined in eq.(\ref{eq:xin}) and (\ref{eq:etan}),
respectively.
\end{theorem}
Remark. Let $\Xi$ be a random variable whose square is subject to a 
$\chi^2$ distribution with freedom 2. 
In accordance with this theorem, $L(X^n)$ converges in distribution to
$L_{\infty}(\Xi)$. 
Hence, 
the rejection region of the most poweful test can be 
asymptotically determined by $L_{\infty}(\Xi)$. 
\begin{proof}
The log density ratio function is given by 
\begin{eqnarray*}
f(X_{i},a,b,c) &=&\log\frac{p(X_i|a,b,c)}{p(X_i|0,0,1)}
\\
&=& \log{\left[1-a + a\sqrt{c}\; e^{g(X_i,b,c)}\right]},
\end{eqnarray*}
where $g(x,b,c)$ is a function defined by 
\[
g(X_i,b,c)=-\left(\frac{c}{2}-\frac{1}{2}\right)X_{i}^{2}+bcX_{i}- \frac{b^{2}c}{2}. 
\]
Hence
the marginal likelihood ratio $L=L(X^n)$ is given by
\begin{eqnarray*} 
L &=&\int_0^1da \int_{D} dbdc\;
\varphi_1(a,b,c) \;\exp\left(\sum_{i=1}^n f(X_i,a,b,c)\right),
\end{eqnarray*}
where
\[
\varphi_1(a,b,c)=\frac{ n}{\sqrt{2} \pi R_0^2}.
\]
Since the integrated region of $(b,c)$ of this integral is $D$, 
$b=O_p(n^{-1/2})$, $c=O_p(n^{-1/2})$. It follows that 
\begin{eqnarray*}
f(X_{i},a,b,c) &=& f|_{(b,c) = (0,1)} + \frac{\partial f}{\partial b}\Big{|}_{(b,c) = (0,1)} b + \frac{\partial f}{\partial c}\Big{|}_{(b,c) = (0,1)}\left(c-1\right) \\
&+& \frac{1}{2}\frac{\partial^{2} f}{\partial b^{2}}\Big{|}_{(b,c) = (0,1)}b^{2} + \frac{1}{2}\frac{\partial^{2} f}{\partial c^{2}}\Big{|}_{(b,c) = (0,1)}\left(c-1\right)^{2}  \\
&+& \frac{\partial^{2} f}{\partial b\partial c}\Big{|}_{(b,c) = (0,1)}b\left(c-1\right) +o_p(1/n),
\end{eqnarray*}
where
\begin{align*}
&\frac{\partial f}{\partial b}\Big{|}_{(b,c) = (0,1)} = aX_{i} \\
&\frac{\partial f}{\partial c}\Big{|}_{(b,c) = (0,1)} = \frac{a}{2}\left(X_{i}^{2} - 1\right)
\\
& \frac{\partial^{2} f}{\partial b^{2}}\Big{|}_{(b,c) = (0,1)} = a\left(X_{i}^{2}-1\right) - a^{2} X_{i}^{2} 
\\
&\frac{\partial^{2} f}{\partial b\partial c}\Big{|}|_{(b,c) = (0,1)} = aX_{i} + \frac{a}{2}\left(1-a\right)X_{i}\left(1-X_{i}^{2}\right) 
\\
&\frac{\partial^{2} f}{\partial c^{2}}\Big{|}_{(b,c) = (0,1)} = -\frac{a}{4}\left(1 + X_{i}^{2}\right) 
- \frac{a}{4}\left(X_{i}^{2}- X_{i}^{4}\right) - \frac{a^{2}}{4}\left(1-X_{i}^{2}\right)^{2}.
\end{align*}
Hence
\begin{eqnarray}
f (X_{i},a,b,c)&=& ab X_{i} + \frac{a(c-1)}{2}\left(X_{i}^{2} - 1\right) + \frac{1}{2}\left\{a\left(X_{i}^{2}-1\right)-a^{2}X_{i}^{2}\right\}b^{2} \nonumber \\
&&+ \left\{aX_{i} + \frac{a}{2}(1-a)X_{i}\left(1-X_{i}^{2}\right)\right\}b\left(c-1\right) \nonumber \\
&&+ \frac{1}{2}\left[-\frac{a}{4}\left(1 + X_{i}^{2}\right) - \frac{a}{4}\left(X_{i}^{2} - X_{i}^{4}\right)
-\frac{a^{2}}{4}\left(1-X_{i}^{2}\right)^{2}\right]\left(c-1\right)^{2}  \nonumber \\
&& +o_p(1/n).
\end{eqnarray}
Note that the order of the quadratic forms of $(b,c-1)$ is $1/n$ and
\begin{eqnarray*}
(1/n)\sum_i X_i&=& o_p(1), \\
(1/n)\sum_i X_i^2&=& 1+o_p(1), \\
(1/n)\sum_i X_i^3&=& o_p(1), \\
(1/n)\sum_i X_i^4&=& 3+o_p(1). 
\end{eqnarray*}
The log likelihood ratio function is given by 
\begin{eqnarray*}
\sum_{i=1}^n f (X_{i},a,b,c)&=& ab \sum_{i=1}^nX_{i} 
+ \frac{a(c-1)}{2}\sum_{i=1}^n\left(X_{i}^{2} - 1\right) 
- \frac{n}{2}a^{2}b^{2} \nonumber \\
&& - \frac{n}{4}
a^{2}\left(c-1\right)^{2}  +o_p(1).
\end{eqnarray*}
Let us define $(r,\theta)$ by
\begin{eqnarray*}
b&=&r\cos\theta,\\
c&=&1+\sqrt{2}r\sin\theta.
\end{eqnarray*}
Then by using eq.(\ref{eq:xin}) and (\ref{eq:etan}), 
\begin{eqnarray*}
\sum_{i=1}^n f (X_{i},a,b,c)&=& -\frac{n}{2}a^2r^2+\sqrt{n} ar(\xi_n\cos\theta+\eta_n\sin\theta)+o_p(1)
\\
&=& -\frac{n}{2}a^2r^2+\sqrt{n}ar\sqrt{\xi_n^2+\eta_n^2}\sin(\theta+\theta_0)+o_p(1),
\end{eqnarray*}
where $\theta_0$ is a random variable which satisfies $\tan \theta_0=\xi_n/\eta_n$. 
By using the notation $R=R_0/\sqrt{n}$, 
the log marginal likelihood ratio can be written as 
\begin{eqnarray*}
L&=& \int_0^1da \int_{D}dbdc\;\frac{n}{\sqrt{2}\pi R_0^2}\exp\left(\sum_i f_i(a,b,c)\right)
\\
&=&\int_0^1 da \int_0^{R}\frac{2r\;dr}{R^2}\int_0^{2\pi}\frac{d\theta}{2\pi}
\exp\left(-\frac{n}{2}a^2r^2+ar\sqrt{\xi_n^2+\eta_n^2}\sin(\theta+\theta_0)+o_p(1)
\right)
\\
&=&\int_0^1 da \int_0^{R}\frac{2r\;dr}{R^2}\int_0^{2\pi}\frac{d\theta}{2\pi}
\exp\left(-\frac{n}{2}a^2r^2+ar\Xi_n\sin(\theta)+o_p(1)
\right).
\end{eqnarray*}
Then by replacing 
$r=\ell/\sqrt{n}$ with $dr=d\ell/\sqrt{n}$, it follows that 
\begin{eqnarray*}
L&=&\int_0^1 da \int_0^{R_0}\frac{2 \ell\;d\ell}{R_0^2}\int_0^{2\pi}\frac{d\theta}{2\pi}
\exp\left(-\frac{1}{2}a^2\ell^2+a\ell\Xi_n\sin(\theta)+o_p(1)
\right). 
\end{eqnarray*}
We define $L_{\infty}(\Xi_n)$ by 
\begin{eqnarray*}
L_{\infty}(\Xi_n)&=&\int_0^1 da \int_0^{R_0}\frac{2 \ell\;d\ell}{R_0^2}\int_0^{2\pi}\frac{d\theta}{2\pi}
\exp\left(-\frac{1}{2}a^2\ell^2+a\ell\Xi_n\sin(\theta)
\right).
\end{eqnarray*}
Then, the convergence in probability
$L(X^n)-L_{\infty}(\Xi_n)\rightarrow 0$ holds. 
$L_{\infty}(\Xi_n)$ can be rewritten as 
\begin{eqnarray*}
L_{\infty}(\Xi_n)=\int_0^1 da \int_0^{R_0}\frac{2 \ell\;d\ell}{R_0^2}\int_0^{\pi}\frac{d\theta}{2\pi}
\exp\left(-\frac{1}{2}a^2\ell^2\right)
\cosh\left(a\ell\Xi_n\sin(\theta)
\right).
\end{eqnarray*}
By using 
\[
t=a^2\ell^2,
\]
the random variable $L_{\infty}(\Xi_n)$ can be also rewritten as
\begin{eqnarray*}
L_{\infty}(\Xi_n)&=&\int_0^1 da \int_0^{a^2R_0^2}\frac{\;dt}{a^2 R_0^2}\int_0^{\pi}\frac{d\theta}{2\pi}
\exp\left(-\frac{t}{2}\right)
\cosh\left(\sqrt{t}\Xi_n\sin(\theta)
\right),
\\
&=& \int_0^{R_0^2} dt \int_{\sqrt{t}/R_0}^1 \frac {da} {a^2 R_0^2}\int_0^{\pi}\frac{d\theta}{2\pi}
\exp\left(-\frac{t}{2}\right)
\cosh\left(\sqrt{t}\Xi_n\sin(\theta)
\right),
\\
&=& \frac{1}{2R_0^2}\int_0^{R_0^2} dt 
\left(
R_0/\sqrt{t}-1
\right)
\exp\left(-\frac{t}{2}\right)
I_0(\sqrt{t}\Xi_n),
\end{eqnarray*}
which completes the theorem. 
\end{proof}

As well as the results we obtained in the previous sections, the asymptote of $L$ does not explicitly depend on the sample size $n$. The reason for this behavior is the same as in the previous cases, the critical scaling $r \propto n^{-1/2}$.

Let us validate the asymptote we derived above. First, we will numerically calculate the behavior of $\Xi$, by using a sample generated from the standard normal distribution. 

From a well-known result on the percentile point of the $\chi^{2}$ result, we obtain the $10\%$ percentile as $\Xi \simeq 2.146$, and the $5\%$ percentile as $\Xi \simeq 2.448$, and the $1\%$ percentile as $\Xi \simeq 3.035$.
Therefore, we can construct a hypothesis test by using $\Xi$ as a test statistics.

Next, we describe our numerically validation of the asymptotic distribution of $L$ when the sample size is finite.

We firstly prepared the 10000 sample sets, whose size is denoted by $n$, and conducted the validation for three different values of $n$, i.e., $n=200,400,800$. 
We calculated the $L$ by using $\Xi$ calculated from the finite sample. Here, we set $R_{0}$ as 1.
Then, we estimated the level numerically in each case and compared them with the levels of case 2.

Table 3  shows the result. Compared with the previous cases that used simpler models, in the present case, the numerically calculated levels slightly deviate from the theoretical values derived from the asymptote.
But we can see that as the $n$ becomes larger, the numerically calculated levels approach the theoretical value, and we can conclude that they match well and the asymptote we derived in the previous section works well even in the case of finite $n$.

\begin{table}[h]
\begin{center}
 \caption{Comparison of levels derived from asymptote and those numerically calculated levels}
\vspace{5mm}

 \begin{tabular}{|c|c|r|r|} \hline
 level & 10\% & 5\% & 1\% \\\hline
 rejection region $L>r$   & r=0.550& r=0.581& r=0.659\\\hline
numerically calculated level(n=100) & 9.53\% & 4.81\%& 1.21\% \\\hline
numerically calculated level(n=200) & 9.72\% & 4.64\%& 0.88\% \\\hline
numerically calculated level(n=400) & 10.04\% & 4.91\%& 0.79\% \\\hline
numerically calculated level(n=800) & 10.33\% & 5.09\%& 1.03\% \\\hline
\end{tabular}

 \end{center}
\end{table}

\subsection{Comments from the perspective of the singular learning theory}

To conclude this section, let us mention the relation between the result obtained above and the general asymptotic form of the log marginal likelihood of the singular model, which is derived from the theory of algebraic geometry\citep{Watanabe2001}.

In Theorem 2, we derived the asymptotic form of $L$ and saw that $L$ did not depend on the sample size $n$ as a result of the scaling law $B \propto n^{-1/2}$ that we applied. 

We can consider another scaling $B \propto n^{-\alpha}$, where $\alpha > 0$ is a constant. As long as $\alpha \leq \frac{1}{2}$, we can calculate the asymptotic form of $L$ in the same way as the derivation of Theorem 2. The result is as follows.

\begin{equation}
L = \frac{1}{B_{0}n^{\frac{1}{2}-\alpha}}\int_{0}^{B_{0}^{2}n^{1-2\alpha}}\frac{1}{\sqrt{t}}\log{\left[\frac{B_{0}^{2}n^{1-2\alpha}}{t}\right]}e^{-t/2}\cosh{\xi \sqrt{t}} dt
\end{equation}

We can immediately obtain the log marginal likelihood ratio $F = -\log L$.

\begin{equation}
F = \left(\frac{1}{2}-\alpha\right) \log n - \left(1-2\alpha\right)\log\left(\log n\right) + o_{p}(\log\left(\log n\right))
\end{equation}

From the general theory, the asymptotic form of log marginal likelihood becomes

\begin{equation}
F = \lambda \log n- \left(m-1\right)\log{\left(\log n\right)} + o_{p}(\log\left(\log n\right))
\end{equation}

We can see that our result corresponds to $\lambda = \frac{1}{2}-\alpha$ and $m = \left(2-2\alpha\right)$. The sample size's dependency on the support of the prior affects the real canonical log threshold $\lambda$ and the multiplicity $m$. In this paper, we treated $\alpha=\frac{1}{2}$  as a ``critical'' case, where the $\lambda$ and $m$ effectively vanish. In such a case, the main term of $F$ becomes stochastic. This is why it can be difficult to apply conventional Bayes factor-based testing to such a case. 

Let us comment more on the scaling $n^{-1/2}$. The Kullback-Leibler divergence $K$ between the null hypothesis and the alternative hypothesis can be easily calculated.

\begin{equation}
K(a.b) = \int p(x|(0,0)) \log{\frac{p(x|(0,0))}{p(x|w)}} dx= - \frac{1}{2}a^{2}b^{2}
\end{equation}

Here, $nK(a,b)$ is nothing other than the leading term of the $H(a.b)$. 

In the proof of the Theorem 2,  we mainly considered that $b \propto n^{-1/2}$, and $a \sim \mathcal{O}(1)$. The meaning of this setup is clear, the center of the mixed distribution deviates from the origin , as much as the variance of the distribution, and the null and alternative hypothesis are hard to discriminate. 
 
As a result of this scaling, both $na^{2}b^{2}$ and $\sqrt{n}ab\xi$ becomes $\mathcal{O}(1)$, and this result in the $n$-independent asymptote of $L$. 
 
However, as we can be easily seen, this ``scaling'' is not unique. So long as $ab \sim n^{-1/2}$ and $bX_{i}-\frac{b^{2}}{2}$ is  small enough that the Taylor expansion of the exponential is valid, a proof similar to the one above can be constructed. For example, a scaling such as $a\sim n^{-1/4}$ and $b \sim n^{-1/4}$ will lead to the same results.

The important point here is that this can be understood as a Taylor expansion around the singularity $ab=0$, and the deviation is described as a power of $ab$, not of $b$. 

As we saw above, in this delicate situation, a hypothesis test based on the stochastic behavior of $L$ works well, and to construct it, we need to find the singularity (in our setting, $ab=0$)and an appropriate scaling (in our setting, $ab \sim n^{-1/2})$is essentially important. 

Therefore, to construct the hypothesis test using singular models, we should keep in mind the effect of the singularity, and consider whether the case under consideration is ``delicate'' or not, by computing the Kullback-Leibler divergence between the null hypothesis and the alternative hypothesis. The scaling is determined by the form of the Kullback-Leibler divergence that consists of a polynomial for each parameters. From the perspective of the singular learning theory, this is nothing other than the relation between the real log canonical threshold (RCLT) $\lambda$ and the representation of the parameters in the model.

\section{Conclusion}
\label{sec:conc}
In this paper, we theoretically studied the test of homogeneity for normal mixtures in terms of the Bayesian  framework, for the first time.

By applying the mathematical technique developed for the analysis of singular models and by appropriately  scaling from the singularity, we derived the asymptotic behavior of the marginal likelihood ratio for several forms of the prior.  These forms are clearly different from the conventional $\chi^{2}$ distribution, as an effect of the singularity in the parameter space, but their stochastic behavior can be described as a function of random variables that obey the normal distributions. We constructed a hypothesis test based on these results and numerically validated their effectiveness.

The merits of our treatment, based on the Bayesian learning theory for singular models, are as follows.

First, the test statistics that we analyzed was the marginal likelihood ratio and as a result of this, the hypothesis test using it is guaranteed to be the most powerful test. Second, compared with other methods using the value of the (log) likelihood ratio, such as Bayes factor based ones, the hypothesis test based on the stochastic behavior of the marginal likelihood ratio is valid even when the null hypothesis and the alternative one are hard to discriminate, as we saw in Section 4. The stochastic behavior of the test statistics we derived can be described as a function of  the probability variables obeying well-known probability distributions. From the practical perspective, this gives us a clear and easy-to-use formalism. 

We should note that the construction of our hypothesis test is not possible until the stochastic behavior of the marginal likelihood ratio is theoretically derived. As far as we know, this is the first time a concrete form was derived.  The results are of the mathematical theory that enable us to treat singular models properly.

To conclude our discussion,  we should note that in Bayesian learning theory, the study of hypothesis tests is not sufficient and there is much that remains to be studied.  We believe that our method is very general, and that it can be applied to various singular models. This direction of study could be of practical value. We also  believe that it is also important to study the methods of approximating the log marginal likelihood ratios with high accuracy. One candidate for this is variational Bayes, which is an efficient way to approximate the posterior distribution. However, the theory of hypothesis test based on variational Bayes is still insufficient. Therefore, in the future, we should study how to apply it to a Bayesian hypothesis test.

\bibliographystyle{unsrtnat}

\bibliography{list}

\begin{thebibliography}{15}
\providecommand{\natexlab}[1]{#1}
\providecommand{\url}[1]{\texttt{#1}}
\expandafter\ifx\csname urlstyle\endcsname\relax
  \providecommand{\doi}[1]{doi: #1}\else
  \providecommand{\doi}{doi: \begingroup \urlstyle{rm}\Url}\fi

\bibitem[Pearson(1894)]{Pearson1894}
K~Pearson.
\newblock Iii. contributions to the mathematical theory of evolution.
\newblock \emph{Philosophical Transactions of the Royal Society of London A:
  Mathematical, Physical and Engineering Sciences}, 185:\penalty0 71--110,
  1894.
\newblock ISSN 0264-3820.
\newblock \doi{10.1098/rsta.1894.0003}.

\bibitem[McLachlan and Peel(2000)]{Mclachlan2000}
G.~J. McLachlan and D.~Peel.
\newblock \emph{Finite mixture models}.
\newblock Wiley Series in Probability and Statistics, New York, 2000.

\bibitem[Hartigan(1985)]{Hartigan1985}
J.~A. Hartigan.
\newblock A failure of likelihood asymptotics for normal mixtures.
\newblock \emph{Proceedings of the Barkeley Conference in Honor of Jerzy Neyman
  and Jack Kiefer, 1985}, 2:\penalty0 807--810, 1985.

\bibitem[Liu and Shao(2003)]{Liu2003}
Xin Liu and Yongzhao Shao.
\newblock Asymptotics for likelihood ratio tests under loss of identifiability.
\newblock \emph{Ann. Statist.}, 31\penalty0 (3):\penalty0 807--832, 06 2003.
\newblock \doi{10.1214/aos/1056562463}.

\bibitem[Garel(2001)]{Garel2001}
Bernard Garel.
\newblock Likelihood ratio test for univariate gaussian mixture.
\newblock \emph{Journal of Statistical Planning and Inference}, 96\penalty0
  (2):\penalty0 325 -- 350, 2001.
\newblock ISSN 0378-3758.

\bibitem[Chen et~al.(2001)Chen, Chen, and Kalbfleisch]{Chen2001}
Hanfeng Chen, Jiahua Chen, and John~D. Kalbfleisch.
\newblock A modified likelihood ratio test for homogeneity in finite mixture
  models.
\newblock \emph{Journal of the Royal Statistical Society. Series B (Statistical
  Methodology)}, 63\penalty0 (1):\penalty0 19--29, 2001.
\newblock ISSN 13697412, 14679868.

\bibitem[Chen et~al.(2004)Chen, Chen, and Kalbfleisch]{Chen2004}
Hanfeng Chen, Jiahua Chen, and John~D. Kalbfleisch.
\newblock Testing for a finite mixture model with two components.
\newblock \emph{Journal of the Royal Statistical Society: Series B (Statistical
  Methodology)}, 66\penalty0 (1):\penalty0 95--115, 2004.
\newblock \doi{10.1111/j.1467-9868.2004.00434.x}.

\bibitem[Chen and Li(2009)]{Chen2009}
Jiahua Chen and Pengfei Li.
\newblock Hypothesis test for normal mixture models: The em approach.
\newblock \emph{The Annals of Statistics}, 37\penalty0 (5A):\penalty0
  2523--2542, 2009.
\newblock ISSN 00905364.

\bibitem[Chen et~al.(2012)Chen, Li, and Fu]{Chen2012}
Jiahua Chen, Pengfei Li, and Yuejiao Fu.
\newblock Inference on the order of a normal mixture.
\newblock \emph{Journal of the American Statistical Association}, 107\penalty0
  (499):\penalty0 1096--1105, 2012.
\newblock \doi{10.1080/01621459.2012.695668}.

\bibitem[Charnigo and Sun(2004)]{Charnigo2004}
Richard Charnigo and Jiayang Sun.
\newblock Testing homogeneity in a mixture distribution via the l2 distance
  between competing models.
\newblock \emph{Journal of the American Statistical Association}, 99\penalty0
  (466):\penalty0 488--498, 2004.

\bibitem[Chauveau et~al.(2017)Chauveau, Garel, and Mercier]{Chauveau2017}
Didier Chauveau, Bernard Garel, and Sabine Mercier.
\newblock {Testing for univariate Gaussian mixture in practice}.
\newblock working paper or preprint, November 2017.

\bibitem[Watanabe(2018)]{Watanabe2018}
Sumio Watanabe.
\newblock \emph{Mathematical Theory of Bayesian Statistics}.
\newblock Chapman and Hall/CRC, New York, 2018.

\bibitem[Watanabe(2010)]{Watanabe2010}
Sumio Watanabe.
\newblock Asymptotic equivalence of bayes cross validation and widely
  applicable information criterion in singular learning theory.
\newblock \emph{J. Mach. Learn. Res.}, 11:\penalty0 3571--3594, December 2010.
\newblock ISSN 1532-4435.

\bibitem[Watanabe and Amari(2003)]{Watanabe2003}
Sumio Watanabe and Shun-ichi Amari.
\newblock Learning coefficients of layered models when the true distribution
  mismatches the singularities.
\newblock \emph{Neural Computation}, 15\penalty0 (5):\penalty0 1013--1033,
  2003.
\newblock \doi{10.1162/089976603765202640}.

\bibitem[Watanabe(2001)]{Watanabe2001}
Sumio Watanabe.
\newblock Algebraic analysis for nonidentifiable learning machines.
\newblock \emph{Neural Computation}, 13\penalty0 (4):\penalty0 899--933, 2001.

\end{thebibliography}
\end{document}